\documentclass{ws-jktr}
\usepackage{wrapfig}

\title{Symmetric Union Presentations for 2-Bridge Ribbon Knots}

\begin{document}

\markboth{Christoph Lamm}
{Symmetric Union Presentations for 2-Bridge Ribbon Knots}

\author{Christoph Lamm}

\address{R\"{u}ckertstr. 3, 65187 Wiesbaden, Germany \\ 
e-mail: christoph.lamm@web.de}

\maketitle

\vspace{2.6cm}

\begin{abstract}
We prove that all 2-bridge ribbon knots are symmetric unions. 
\end{abstract}

\section{Overview and definitions}
Symmetric unions have been defined as generalisations of Kinoshita-Terasaka's construction in 1957
(\cite{KinoshitaTerasaka}, \cite{Lamm}).
They are given by diagrams which look like the connected sum of a knot and its mirror image with 
additional twist tangles inserted near the symmetry axis. Because all symmetric unions are ribbon knots, 
we can ask how big a subfamily of ribbon knots they form. It is, for instance, known that all 21 
(non-trivial, prime) ribbon knots with crossing numbers $\le 10$ are symmetric unions. 
In this article we prove that all 2-bridge ribbon knots are symmetric unions.

\begin{figure}[hbtp]
\centering
\includegraphics[scale=1.1]{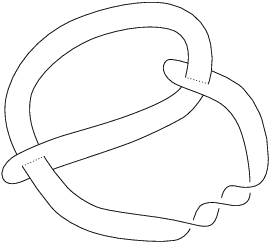}\hspace{2cm}
\includegraphics[scale=1.1]{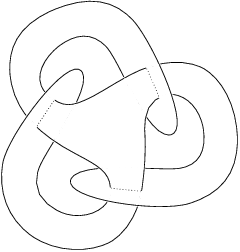}
\caption{Examples of disks with ribbon singularities}
\label{ribbon}
\end{figure}

This article is based on a talk at the Joint Meeting of AMS and DMV in June 2005 at Mainz University (Germany).
We include a postscriptum containing comments and later developments. Note that in 2005 it
was not yet known that the three Casson-Gordon families of 2-bridge ribbon knots  are a complete 
list of all 2-bridge ribbon knots. This was proved by Paolo Lisca in 2006.

\begin{definition}
A knot $K$ in $S^3$ is a {\sl ribbon knot} if it bounds an immersed disk in $S^3$ with only ribbon singularities.
\end{definition}

Recall that $K$ is a {\sl slice knot} if it bounds a smoothly embedded 2-disk $D^2$ in $B^4$. 
Every ribbon knot is slice. A notorious question in knot theory is whether the converse is true.

\begin{figure}[hbtp]
\centering
\includegraphics[scale=1.0]{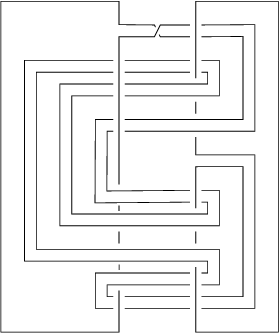}
\caption{\small A ribbon presentation for the knot $10_{87}$}
\label{kawa}
\end{figure}

There are 21 non-trivial prime ribbon knots with crossing numbers $\le$ 10. A list is shown in 
Appendix F.5 in the book  "A survey of knot theory" \cite{Ka}. An example can be seen in Figure \ref{kawa}.

\begin{figure}[hbtp]
\centering
\includegraphics[scale=1.25]{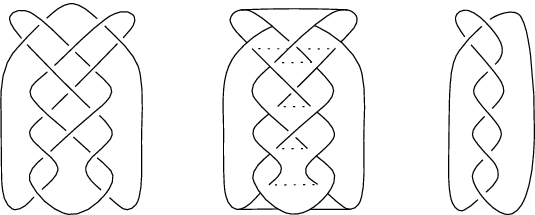}
\caption{\small The knot $10_{87}$ as a knot diagram (left), as a symmetric ribbon with twists (center)
and the partial knot diagram (right)}
\label{10_87}
\end{figure}

Symmetric unions yield another way of showing that knots are ribbon.
As a general reference we refer to \cite{Lamm}. Before we give the definition of symmetric unions we mention 
that all 21 non-trivial prime ribbon knots with crossing number less or equal to 10 are symmetric unions.
As an example we show in Figure \ref{10_87} a symmetric union for $10_{87}$ a case which was missing 
from my list in 1998 \cite{Lamm}.

\begin{figure}[hbtp]
\centering
\includegraphics[scale=0.7]{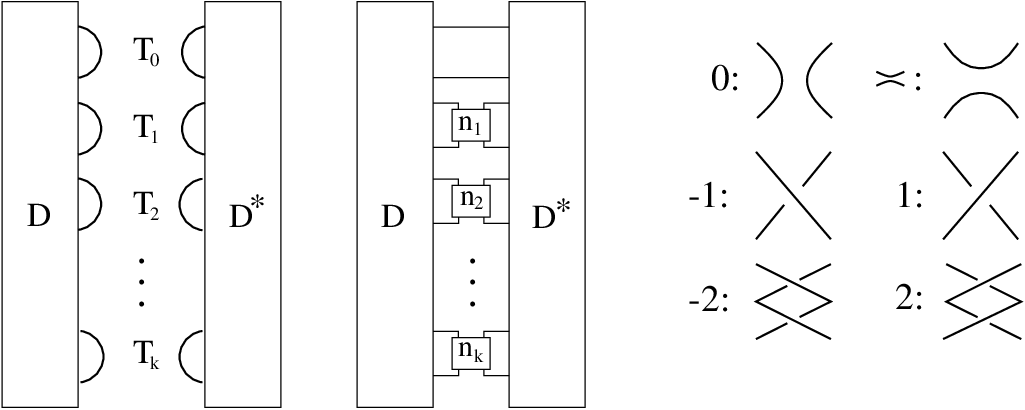}
\caption{Definition of symmetric unions}
\label{sudef}
\end{figure}

\begin{definition}
Let $D$ be an unoriented knot diagram and $D^*$ the diagram $D$ reflected at an axis in the plane.
If, as in Figure \ref{sudef}, we insert the tangle $\asymp$ and twist tangles $n_i$ we call the result a
{\sl symmetric union} of $D$ and $D^*$ with twist parameters $n_i$, $i = 1, \ldots , k$. 
The {\sl partial knot} $\hat K$ of a symmetric union of $D$ and $D^*$ is the knot given by the diagram $D$. 
\end{definition}

\bigskip\noindent
{\bf Examples:}

\noindent
1.) As an example we consider again the diagram in Figure \ref{10_87}. 
Here the partial knot is the knot $6_1$, shown in the right of the figure.

\noindent
2.) If all $n_i=0$, then we get the well-known symmetric ribbon for $\hat K \sharp -\hat K$ (we use the
notation $-K$ for the mirror image of $K$, with reversed orientation).

\bigskip 

\bigskip
\section{Properties of symmetric unions}
The following Theorem \ref{matrix} was already known to Kinoshita and Terasaka 
in the case that only one twist tangle is inserted, see \cite{KinoshitaTerasaka}. 
The theorem is proved for instance by using the matrix definition of the Alexander 
polynomial and the Goeritz matrix \cite{Lamm}.

\begin{theorem}\label{matrix}
The Alexander polynomial of a symmetric union of $D$ and $D^*$ 
depends only on the parity of the twist numbers $n_i$.
The determinant of a symmetric union of $D$ and $D^*$ is independent of the twist 
numbers and equals the square of the determinant of the partial knot.
\end{theorem}

\bigskip
In this article we are interested mainly in the ribbon property of symmetric unions:

\begin{theorem}\label{symribbon}
Symmetric unions are ribbon knots.
\end{theorem}

\begin{figure}[hbtp]
\centering
\includegraphics[scale=1.0]{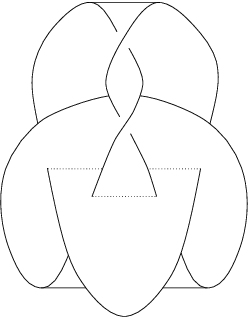} \hspace{2.5cm}
\includegraphics[scale=0.09]{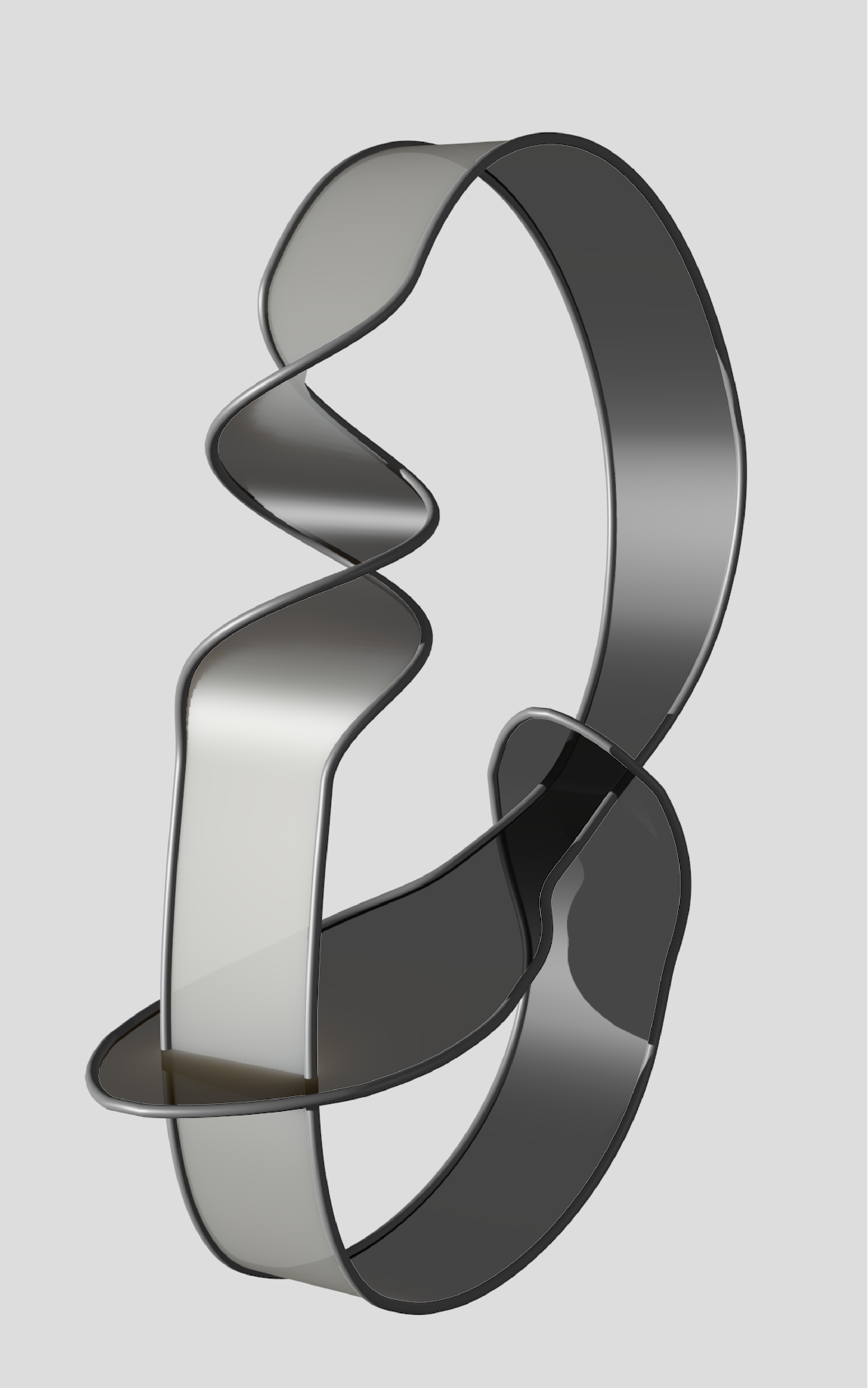}
\caption{Symmetric unions are ribbon knots (left: a diagram of $8_{20}$ showing the symmetric ribbon,
right: a 3D rendering)}
\label{bandbeisp}
\end{figure}

The proof uses the same construction as for $K \sharp -K$ with additional twists in the ribbon, 
see Figure \ref{bandbeisp} for an example. The articles \cite{Lamm}, \cite{EisermannLamm2007} 
and \cite{EisermannLamm2011} contain the proof and more details on symmetric unions.

\begin{figure}[hbtp]
\centering
\includegraphics[scale=0.65]{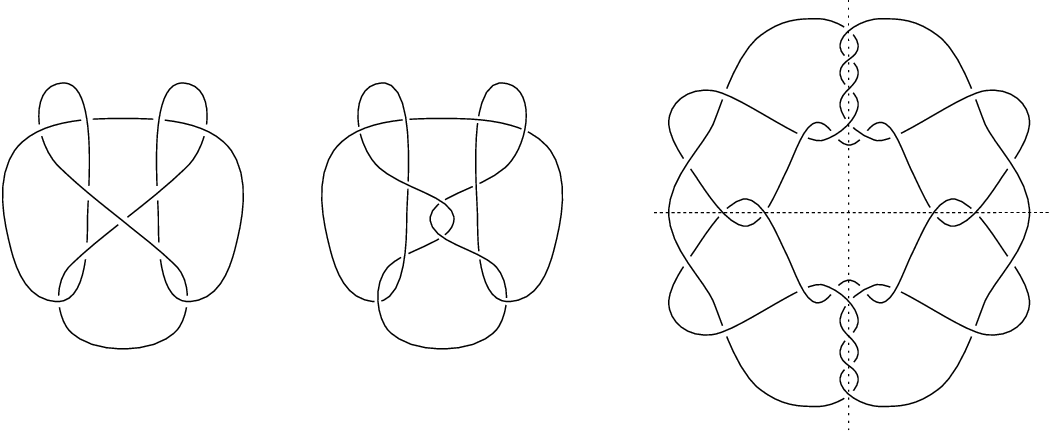}
\caption{The knot $10_{153}$ (left), the Kinoshita-Terasaka knot (center) and 
a symmetric union with two symmetry axes used by Kanenobu (right)}
\label{det1}
\end{figure}

\noindent
Further examples:
\begin{itemize}
\item The knot $10_{153}$ and the Kinoshita-Terasaka knot are symmetric unions of the trivial knot, see Figure \ref{det1}. 
Hence they have determinant 1. For the Kinoshita-Terasaka knot the twist parameter is an even number,
hence its Alexander polynomial is equal to 1.

\item In 2004 Taizo Kanenobu used the chiral knot in the right of Figure \ref{det1} in order to present a knot 
whose chirality is not detected by the Links/Gould invariant, see \cite{deWitLinks}, \cite{IshiiKanenobu}.
This knot is a symmetric union with respect to two different axes, resulting in different partial knots.
A similar diagram can be found in Example 3.1 in \cite{EisermannLamm2007}.
\end{itemize}

\section{2-bridge knots: definition and the Casson-Gordon families}
Recall the definition of the bridge number of a knot: if $v \in \mathbb{S}^2$ is a unit vector in $\mathbb{R}^3$ 
and $K$ is a knot, then let $b_v(K)$ be the number of maxima of the orthogonal projection of $K$ on the line 
spanned by $v$. Then the {\sl bridge number} of $K$ is:
$$b(K):=\min_{K' \sim K} \min_v b_v(K')$$
Because knots with bridge number 1 are trivial, the simplest cases occur for bridge number 2. 
They were studied by:
\begin{itemlist}
\item Bankwitz/Schumann (Viergeflechte, 1934),
\item Schubert (Knoten mit 2 Br\"ucken, 1956) and
\item Conway (rational knots, 1970).
\end{itemlist}

We mention some properties of 2-bridge knots:
Their double branched coverings are lens spaces $L(p,q)$ with $p$ equal to the determinant of the knot.
In the plait normal form $C(a_1,\ldots,a_n)$, also called {\sl Conway notation}, 
the numbers of half-twists $a_i$ are related to the parameters $p$ and $q$ 
by the regular (positive) continued fraction $[a_1, \ldots, a_n] = \frac{p}{q}$ (see \cite{Ka}, Chapter 2.1).
We sometimes call $\frac{p}{q}$ the `fraction' of $K$. It determines the knot type.

The following theorem gives a necessary condition for 2-bridge ribbon knots.

\begin{theorem}[Casson-Gordon, 1974, \cite{CG1}]\label{CassonGordonTheorem}\\
Let $K$ be a 2-bridge knot with fraction $p^2/q$. Denote by $\Delta(a,b)$, where $a, b \in \mathbb{R}$,
the triangle with vertices $(0,0)$, $(a,0)$ and $(a,b)$.
If $K$ is ribbon then we have: 
$$4 (\textnormal{Area} \Delta (pr,\frac{qr}{p}) - \textnormal{Int} \Delta (pr,\frac{qr}{p})) =  
\pm 1, \quad \forall r \in \{1,\ldots,p-1\},$$
where $\textnormal{Int} \Delta(a,b)$ is a weighted count of the integral points lying in $\Delta(a,b)$
(see below for the precise definition).
\end{theorem}

The value of $\textnormal{Int} \Delta (a,b)$ is computed by counting lattice 
points (the set $\mathbb{Z}^2 \cap \Delta (a,b)$) similar to Pick's formula \cite{EisermannLamm2009}: 
interior points count as 1, boundary points as 1/2 and vertices different from (0,0) as 1/4.
The vertex (0,0) is not counted.

\begin{figure}[hbtp]
\centering
\includegraphics[scale=0.55]{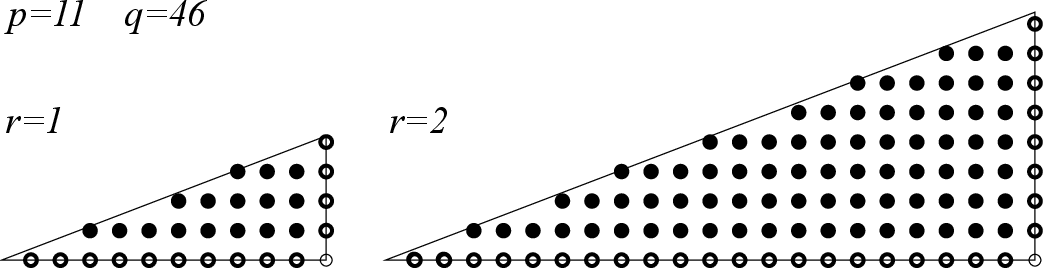}
\caption{Computation of $\textnormal{Int} \Delta (pr,\frac{qr}{p})$}
\label{gordon}
\end{figure}

\noindent
In the example of Figure \ref{gordon}, with $p=11$ and $q=46$,

- for $r=1$ the area of the triangle is 23 and the value of 
$\textnormal{Int} \Delta (pr,\frac{qr}{p}))$ is $23+\frac{1}{4}$, 

- for $r=2$ these values are 92 and $91+\frac{3}{4}$. 

It can be checked that for this knot (with fraction $121/46$) the difference between the area and
the lattice point count is $\pm \frac{1}{4}$ for all $r \in \{1,\ldots,10\}$. It therefore
satisfies the necessary condition.

More generally, it was found that the following three families of 2-bridge knots satisfy the 
necessary condition to be ribbon given in Theorem \ref{CassonGordonTheorem}.

\begin{itemlist}
\item Family 0: $C(a,b,\ldots,w,x,x+2,w,\ldots,b,a)$, with parameters $> 0$,
\item Family 1: $C(2a,2,2b,-2,-2a,2b)$, with $a,b\not=0$,
\item Family 2: $C(2a,2,2b,2a,2,2b)$, with $a,b\not=0$,
\end{itemlist}

\noindent
These three families were for a long time the unique known families of 2-bridge ribbon knots satisfying 
the necessary condition (see \cite{CG1} and \cite{Siebenmann}) and families 0 and 2 where shown to be ribbon
(the status of family 1 remained unclear, see end of Section 4).
It was a conjecture formulated in 1974 that there are no other 2-bridge ribbon (or slice) knots. 
Paolo Lisca proved this in 2006, see \cite{Lisca}.

\section{Symmetric union presentations for the 3 families}
Our main theorem is:

\begin{theorem}\label{maintheorem}
Every knot contained in one of the three families is a symmetric union.
\end{theorem}

\begin{proof}
The proof is given by the following knot diagram transformations. For family 0 we have 
$C(a,b,\ldots,w,x+1,x-1,w,\ldots,b,a) = C(a,b,\ldots,w,x,1,-x,-w,\ldots,-b,-a)$
which is a symmetric union, see Figure \ref{family0}. 
This family was already considered by Kanenobu in 1986 to construct 
knots with the same Jones or Homflypt polynomial, see \cite{Kanenobu}.

\begin{figure}[hbtp]
\centering
\includegraphics[scale=0.6]{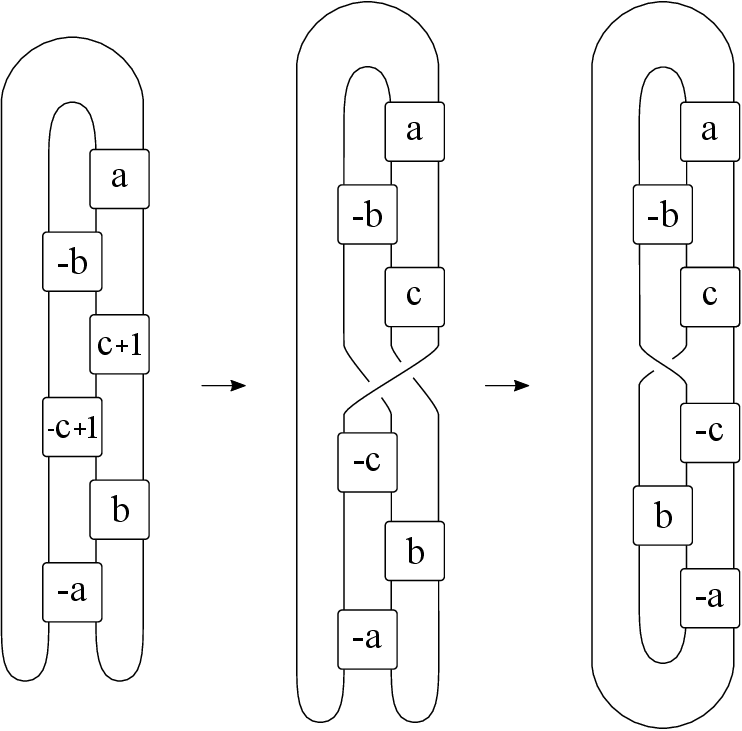}
\caption{Proof for family 0}
\label{family0}
\end{figure}

\bigskip
For families 1 and 2 we have the following diagrams:

\begin{figure}[ht]
\centering
\includegraphics[scale=0.6]{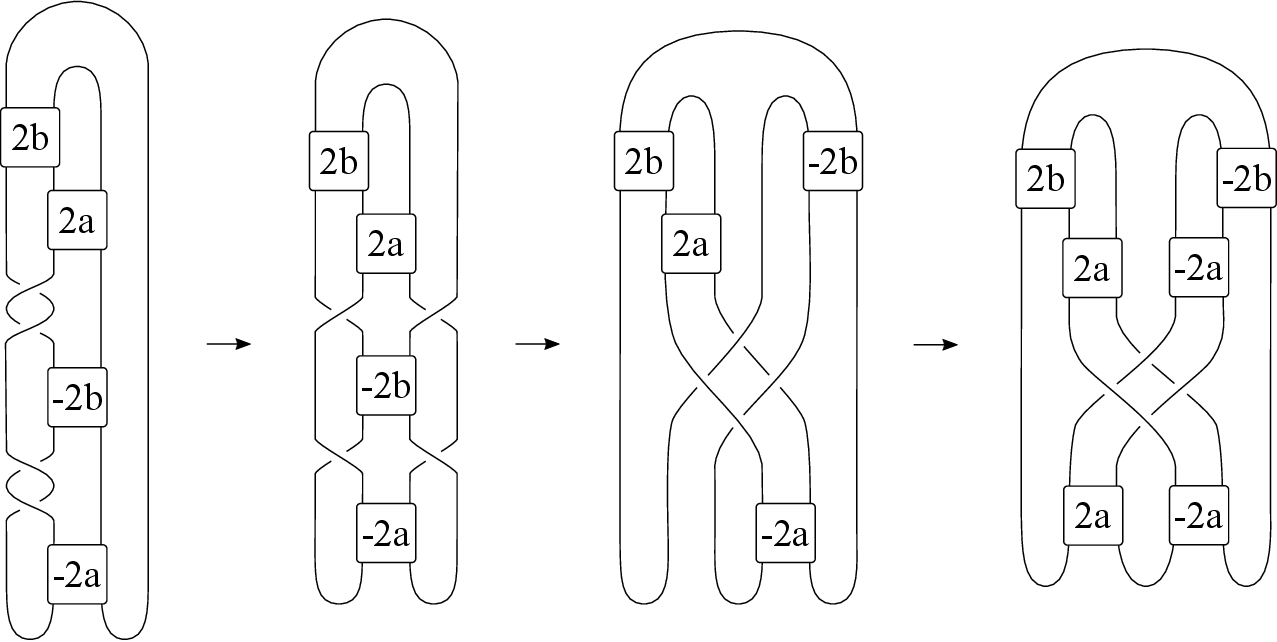}
\caption{Proof for family 1}
\label{family1}
\end{figure}

\clearpage

\begin{figure}[hbtp]
\centering
\includegraphics[scale=0.6]{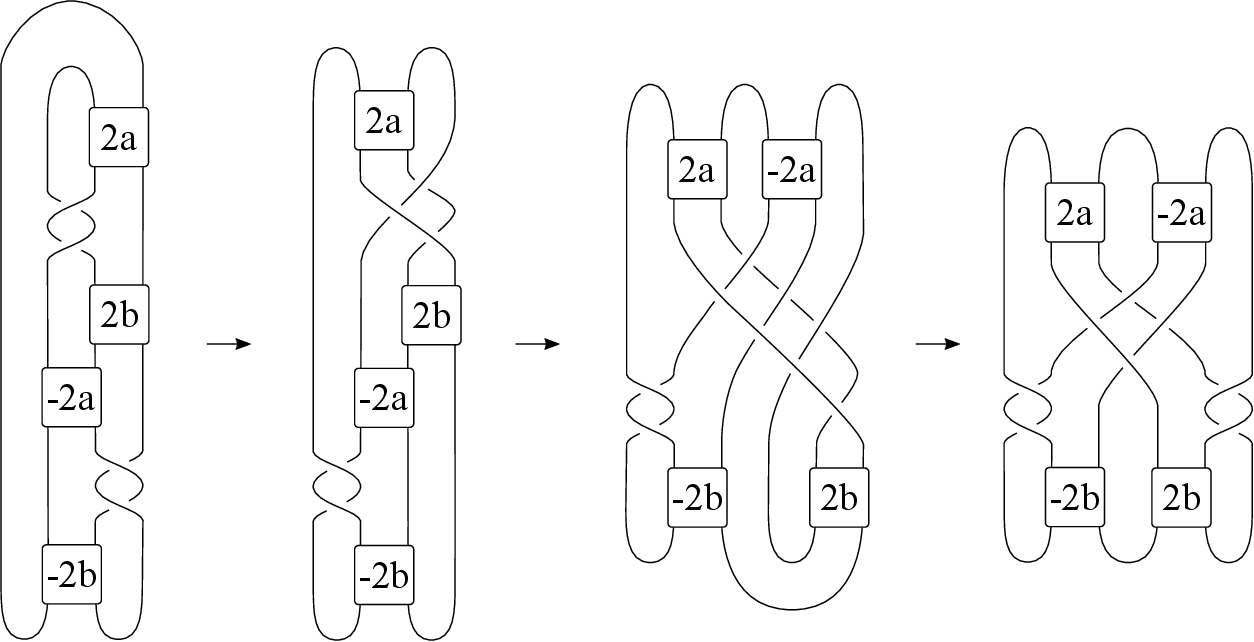}
\caption{Proof for family 2}
\label{family2}
\end{figure}

\end{proof}

Surprisingly, it is necessary to insert extra 2a/-2a twists to obtain the symmetric diagram for family 1.

By Theorem \ref{symribbon} this shows that the knots in the three families are ribbon knots.
For family 1 this seemed, in 2005, to be the first time that their ribbon property was explicitly given.
See also the remark in \cite{Lisca} at the end of the first section and section 3 in \cite{EisermannLamm2009}.

We explain the steps leading to the proof of Theorem \ref{maintheorem}: 
We generated symmetric union diagrams and checked if they represent 2-bridge knots using 
Knotscape \cite{Knotscape} and the list in \cite{deWit}. Four families of diagrams could be
distinguished in this list of symmetric 2-bridge knot diagrams.  For each family we found
a diagram transformation from the symmetric diagram to a plait normal form (corresponding
to the direction from right to left in the above figures) and it was possible to reduce the number of
families to three. The proof reverses the direction of the diagram transformations.

\section{Open questions and a project}

The article of 2005 ended with the following conjectures and a suggested project.

\begin{conjecture}\label{conj1}
Every 2-bridge ribbon knot is contained in one of the three families of Casson and Gordon.

\end{conjecture}

\begin{conjecture}\label{conj2}
Every ribbon knot is a symmetric union.
\end{conjecture}

As mentioned, the first conjecture was proved in \cite{Lisca}. The second is open.

\bigskip\noindent
{\bf Project:}
\begin{itemize}
\item All ribbon knots and symmetric unions with minimal crossing number $\le 10$ are known.
Extend this list to knots with minimal crossing number 11 (is the Conway knot ribbon?).
\end{itemize}

The project was tackled by Axel Seeliger \cite{Seeliger}. He considered knots with crossing number $\le$ 14.
His work and the search for ribbon knots not representable as symmetric unions is  described 
in \cite{Lamm2017} -- taking this information into account we believe that Conjecture \ref{conj2} is wrong. 
Recently, Lisa Piccirillo proved that the Conway knot is not slice \cite{Piccirillo}.

\section{Postscriptum}
As this article is originally from 2005 we comment on later developments.

A detailed article on the (still open) number-theoretic question related to 2-bridge 
ribbon knots is \cite{EisermannLamm2009}. This question asks whether the condition 
in Theorem \ref{CassonGordonTheorem} is sufficient for 2-bridge knots to be ribbon.
Additional related topics are discussed in the recent articles \cite{FellerMcCoy} and \cite{Miller}.

We mention a few more articles on the topic of symmetric unions: by Paolo Aceto \cite{Aceto},
Carlo Collari and Paolo Lisca \cite{CollariLisca}, \cite{CollariLisca2}, Toshifumo Tanaka
\cite{Tanaka2015}, \cite{Tanaka2019} and Feride Ceren Kose \cite{Kose}. Maggie Miller used 
symmetric unions as examples (in particular the diagram for family 2) in Section 6 of \cite{MaggieMiller}.

\medskip \noindent
We end this article with two new projects:

\begin{enumerate}
\item We  propose to study the following relationship: 
two knots are related if they are partial knots for the same knot, given as (two different) symmetric unions. 
Of course, by Theorem \ref{matrix}, two related knots share the same determinant. Are there further restrictions? 
Prove that there are knots with the same determinant which are not related.
\item The second project proposes to find bounds for the crossing number of a symmetric union in terms of 
the crossing number of its partial knot. In the upper direction, it is easy to see that there is no such bound.
For the lower bound, we look at some examples first: 

Figure \ref{10_87} shows a symmetric union for $10_{87}$ with partial knot $6_1$.
Because the crossing number of $6_1 \sharp 6_1^*$ is 12 and three additional crossings 
are added it is remarkable that the crossing number of the resulting knot is as small as 10. 
Motivated by this example we define a new knot invariant: 
for each knot $K$ let $\delta(K)$ be the difference between $2 \cdot c(K)$ and 
the lowest crossing number of a symmetric union with partial knot $K$.
Because  $c(K \sharp -K) \le 2 \cdot c(K)$, this invariant is non-negative.

Examples: a) $\delta(3_1) = 0$, because there are no knots with determinant 9 and crossing number less than 6.
b) $\delta(5_1) = 2$, because $8_8$ and $4_1 \sharp 4_1$ are symmetric unions with partial 
knot $5_1$ and there is no knot with determinant 25 and crossing number less than 8.
c) $\delta(6_1) = 2$, because the smallest crossing number of knots with determinant 81 is 10 and 
we indeed found the case of $10_{87}$ (see Figure \ref{10_87}).

The first unknown value in the table of knots seems to be for $6_3$: 
there are no knots with determinant 169 and crossing number less than 11 and according to \cite{Lamm2017} 
the two possibilities with 11 crossings ($11a164$ and $11a326$) have not yet occured as symmetric unions with 
partial knot $6_3$, hence $\delta(6_3)$ is 0 or 1. 

Taken over all knots, this invariant is not bounded. To see this, use the family from Section 6.4 in \cite{EisermannLamm2011}.
As in general $\delta$ seems to increase for knots with increasing crossing number, we conjecture that there are 
only finitely many prime knots with $\delta(K)=0$, and even stronger that for each $d \ge 0$ there are only
finitely many prime knots with $\delta(K)=d$.
\end{enumerate}

\medskip \noindent
Acknowledgement: 
I thank the referees for their comments which helped to improve the content and the form of this article.


\begin{thebibliography}{99}


\bibitem{CG1}
  A.\,Casson and C.\,Gordon: \textit{Cobordism of classical knots}, in: \`A la Recherche de la Topologie Perdue, 
	ed. by Lucien Guillou and Alexis Marin, Progress in Mathematics, Volume {\bf 62} (1986), 181--199

\bibitem{deWit}
  D.\,De Wit: \textit{The 2-bridge knots of up to 16 crossings},	J. Knot Theory Ramifications {\bf 16} (2007), 997

\bibitem{deWitLinks}
  D.\,De Wit and J.\,Links: \textit{Where the Links-Gould invariant first fails to distinguish nonmutant prime knots},
	J. Knot Theory Ramifications {\bf 16} (2007), 1021

\bibitem{IshiiKanenobu}
  A.\,Ishii and T.\, Kanenobu: 
	\textit{Different links with the same Links-Gould invariant}, Osaka J. Math. {\bf 42} (2005), 273--290

\bibitem{Ka}
  A.\,Kawauchi: A Survey of Knot Theory, Birkh\"auser, 1996.

\bibitem{Kanenobu} T.\,Kanenobu: 
  \textit{Examples of polynomial invariants of knots and links}, Math. Ann. {\bf 275} (1986), 555--572
	
\bibitem{KinoshitaTerasaka}
  S.\,Kinoshita and H.\,Terasaka: \textit{On unions of knots}, Osaka Math. J. {\bf 9} (1957), 131--153
	
\bibitem{Knotscape}
  J.\,Hoste and M.\,Thistlethwaite: \textit{Knotscape -- providing convenient access to tables of knots},
	http://www.math.utk.edu/$\sim$morwen/knotscape.html (1999)
	
\bibitem{Lamm}
  C.\,Lamm: \textit{Symmetric unions and ribbon knots}, Osaka J. Math. {\bf 37} (2000), 537--550
	
\bibitem{Siebenmann}
  L.\,Siebenmann: \textit{Exercices sur les n{\oe}uds rationnels}, unpublished notes (1975)

			
\bigskip
\noindent References added in 2020

\medskip

\bibitem{Aceto}
  P.\,Aceto: \textit{Symmetric ribbon disks}, J. Knot Theory Ramifications {\bf 23} (2014), 1450048

\bibitem{CollariLisca}
  C.\,Collari and P.\,Lisca: \textit{Symmetric union diagrams and refined spin models}, 
	 Canadian Math. Bull., published online in November 2018

\bibitem{CollariLisca2}
  C.\,Collari and P.\,Lisca: \textit{On symmetric equivalence of symmetric union diagrams}, 
	 arXiv: math.GT/1901.10270 (2019)
	
\bibitem{EisermannLamm2007}
  M.\,Eisermann and C.\,Lamm: \textit{Equivalence of symmetric union diagrams}, 
	J. Knot Theory Ramifications {\bf 16} (2007), 879--898
	
\bibitem{EisermannLamm2009}
  M.\,Eisermann and C.\,Lamm: \textit{For which triangles is Pick's formula almost correct?}, 
	Experiment. Math. {\bf 18} (2009), 187--191
	
\bibitem{EisermannLamm2011}
  M.\,Eisermann and C.\,Lamm: \textit{A refined Jones polynomial for symmetric unions}, 
	Osaka J. Math. {\bf 48} (2011), 333--370

\bibitem{FellerMcCoy}
  P.\,Feller and D.\,McCoy: \textit{On 2-bridge knots with differing smooth and topological slice genera},
  Proc. Amer. Math. Soc. {\bf 144} (2016), 5435--5442 

\bibitem{Kose}
  F.C.\,Kose: 
	\textit{A short proof of Tanaka's theorem on composite knots with symmetric union presentations}, arXiv:math.GT/2007.07363 (2020)
	
\bibitem{Lamm2017}
  C.\,Lamm: \textit{The search for nonsymmetric ribbon knots}, Experiment. Math., published online in March 2019

\bibitem{Lisca}
  P.\,Lisca: \textit{Lens spaces, rational balls and the ribbon conjecture}, Geom. Topol. {\bf 11} (2007) 429--472

\bibitem{Miller}
  A.\,Miller: \textit{A note on the topological sliceness of some 2-bridge knots}, 
	Math. Proc. Camb. Phil. Soc. {\bf 164} (2018), 185--191 

\bibitem{MaggieMiller}  
	M.\,Miller: \textit{Extending fibrations on knot complements to ribbon disk complements},
	arXiv:math.GT/1811.09639 (2018)

\bibitem{Piccirillo}
  L.\,Piccirillo: \textit{The Conway knot is not slice}, arXiv:math.GT/1808.02923 (2018)
	
\bibitem{Seeliger}
  A.\,Seeliger: \textit{Symmetrische Vereinigungen als Darstellungen von Bandknoten bis 14 Kreuzungen
	(Symmetric union presentations for ribbon knots up to 14 crossings)}, Diploma thesis, Stuttgart University (2014)

\bibitem{Tanaka2015}
  T.\,Tanaka: \textit{The Jones polynomial of knots with symmetric union presentations}, 
	  J. Korean Math. Soc. {\bf 52} (2015), 389--402
				
\bibitem{Tanaka2019}
  T.\,Tanaka: \textit{On composite knots with symmetric union presentations}, J. Knot Theory Ramifications,
	published online in July 2019
			

\end{thebibliography}
\end{document}